\documentclass[10pt]{amsart}
\usepackage{amsmath}
\usepackage{amscd}
\usepackage{amssymb}
\usepackage{amsthm}
\usepackage{filecontents}
\usepackage[numbers]{natbib}  
\usepackage[draft]{hyperref}

\newenvironment{psmallmatrix}
  {\left(\begin{smallmatrix}}
  {\end{smallmatrix}\right)}

\theoremstyle{plain}
\newtheorem{prop}{Proposition}
\newtheorem{thrm}[prop]{Theorem}

\theoremstyle{definition}

\newtheorem{ex}[prop]{Example}

\makeatletter
\newcommand{\addresseshere}{%
  \enddoc@text\let\enddoc@text\relax
}
\makeatother

\title{Borcherds products of half-integral weight}
\author{Haowu Wang and Brandon Williams}


\address{Max-Planck-Institut f\"ur Mathematik \\ Vivatsgasse 7 \\ 53111 Bonn, Germany}

\email{haowu.wangmath@gmail.com}

\address{Fachbereich Mathematik \\ Technische Universit\"at Darmstadt \\ 64289 Darmstadt, Germany}

\email{bwilliams@mathematik.tu-darmstadt.de}

\begin{document}

\begin{abstract}  We give a necessary and sufficient criterion for the existence of Borcherds products of half-integral weight associated to even lattices that split two hyperbolic planes. In particular we prove that half-integral weight Borcherds products exist for lattices of arbitrary rank.
\end{abstract}

\maketitle

\section{Introduction}

Let $L$ be an even lattice of signature $(l, 2)$. The Borcherds lift \cite{Borcherds} is a multiplicative lifting that takes vector-valued modular forms on $\mathrm{SL}_2(\mathbb{Z})$ with poles at cusps for a finite Weil representation attached to $L$ as input and yields automorphic forms for the orthogonal group of $L$. The \emph{Borcherds products} constructed in this way have many interesting properties. For example, in a neighborhood of any cusp, they are represented by a converging infinite product in which the exponents are Fourier coefficients of their input function. Moreover, the zeros and poles of any Borcherds product lie on rational quadratic divisors and their orders can be read immediately off of the Fourier expansion of the input function. Of particular importance are Borcherds products of \emph{singular weight} $l / 2 - 1$, as the Fourier expansion of any singular-weight modular form is supported only on vectors of norm zero, implying massive cancellation when the product is expanded and suggesting a deeper underlying structure. Indeed, many of the known singular-weight products arise as Weyl denominators of generalized Kac--Moody algebras (see for example \cite{Scheithauer}). 

There are no known examples of singular-weight Borcherds products on lattices $L$ as above with odd $l \ge 5$. More generally, to the authors' knowledge, no examples of Borcherds products of any half-integral weight appear in the literature for lattices with $l \ge 4$. 

In this note we will show that for every $l \in \mathbb{N}$ there exist even lattices of signature $(l, 2)$ that admit holomorphic Borcherds products of half-integral weight. Among even lattices that split two hyperbolic planes, we prove that lattices that admit Borcherds products of half-integral weight have a very simple characterization. (Recall that a lattice $L$ splits two hyperbolic planes if it can be written in the form $L = K \oplus 2U$, where $K$ is an even lattice and where $U$ is unimodular of signature $(1, 1)$, and $2U = U \oplus U$.)

\begin{thrm} Let $L$ be an even lattice of the form $L = K \oplus 2U$ where $K$ is positive-definite and $U$ is unimodular of signature $(1, 1)$ and let $\langle -, - \rangle : L \times L \rightarrow \mathbb{Z}$ be its inner product. The following are equivalent: \\ (i) $L$ admits holomorphic Borcherds products of half-integral weight; \\ (ii) $\langle x,y \rangle \in 8 \mathbb{Z}$ for all $x,y \in K$.
\end{thrm}

\textbf{Acknowledgments.} The authors thank Jan Bruinier and Eberhard Freitag for helpful discussions. H. Wang is grateful to the Max Planck Institute for Mathematics in Bonn for its hospitality and financial support.  B. Williams is supported by a fellowship of the LOEWE focus group Uniformized Structures in Arithmetic and Geometry. 

\section{Preliminaries}

For lattices that split two hyperbolic planes, the Borcherds lift can be conveniently expressed in terms of Jacobi forms. Suppose $K$ is an even positive-definite lattice with dual lattice $K'$. A \emph{weakly holomorphic Jacobi form} of weight $k \in \mathbb{Z}$ and lattice index $K$ is a holomorphic function $\phi : \mathbb{H} \times (K \otimes \mathbb{C}) \rightarrow \mathbb{C}$ satisfying $$\phi \Big( \frac{a \tau + b}{c \tau + d}, \frac{\mathfrak{z}}{c \tau + d} \Big) = (c \tau + d)^k \mathbf{e}\Big( \frac{c}{c \tau + d} Q(\mathfrak{z})\Big) \phi(\tau, \mathfrak{z})$$ and $$\phi(\tau, \mathfrak{z} + \lambda \tau + \mu) = \mathbf{e}(-Q(\lambda)\tau - \langle \lambda, \mathfrak{z} \rangle) \phi(\tau, \mathfrak{z})$$ for all $M = \begin{psmallmatrix} a & b \\ c & d \end{psmallmatrix} \in \mathrm{SL}_2(\mathbb{Z})$ and $\lambda, \mu \in K$. Here $\mathbb{H}$ is the upper half-plane and $\mathbf{e}(x) = e^{2\pi i x}$. The qualifier \emph{weakly holomorphic} means that the Fourier expansion may be a Laurent series in the variable $q$: $$\phi(\tau, \mathfrak{z}) = \sum_{n \gg  -\infty} \sum_{\ell \in K'} c(n, \ell) q^n \zeta^{\ell}, \;  \; \zeta^{\ell} := e^{2\pi i \langle \ell, \mathfrak{z} \rangle}.$$ We do not impose any condition on $\ell$; however the transformation laws imply that for any fixed $n$ only finitely many terms appear in the sum $\sum_{\ell \in K'} c(n,r) \zeta^{\ell}$. 

Now suppose $\phi$ has weight $0$ and integral Fourier coefficients. The \emph{Borcherds lift} of $\phi$ is a meromorphic modular form of weight $c(0, 0) / 2$ on the Type IV domain attached to the lattice $L = K \oplus 2U$. In terms of the tube domain model $$\mathbb{H}_K = \{Z = (\tau, \mathfrak{z}, w) \in \mathbb{H} \times (K \otimes \mathbb{C}) \times \mathbb{H}: \; Q(\mathrm{im}(\mathfrak{z})) < (\mathrm{im} \, \tau) \cdot (\mathrm{im}\, w) \},$$ writing $q = e^{2\pi i \tau}$, $s = e^{2\pi i w}$, and $r^{\ell} = e^{2\pi i \langle \ell, \mathfrak{z} \rangle}$, the Borcherds lift of $\phi$ is represented locally as the product $$q^A r^B s^C \prod_{(n, \ell, m) > 0} (1 - q^n r^{\ell} s^m)^{c(nm, \ell)},$$ where $(n, \ell, m) > 0$ is a positivity condition with respect to a \emph{Weyl chamber} and $(A, B, C)$ is the associated \emph{Weyl vector}. (See Theorem 4.2 of \cite{Gritsenko} for details.) 

This formulation and Borcherds' vector-valued modular forms are related by the theta decomposition. If we write $$\phi(\tau, \mathfrak{z}) = \sum_{\gamma \in K' / K} f_{\gamma}(\tau) \Theta_{K, \gamma}(\tau, \mathfrak{z}), \;  \text{where} \; \Theta_{K, \gamma}(\tau, \mathfrak{z}) = \sum_{\ell \in K + \gamma} q^{Q(\ell)} \zeta^{\ell},$$ then the association $$\phi \leftrightarrow F(\tau) = \sum_{\gamma \in K' / K} f_{\gamma}(\tau) \mathfrak{e}_{\gamma}$$ (where $\mathfrak{e}_{\gamma}$ are the basis elements of $\mathbb{C}[K'/K]$) defines an isomorphism between weakly holomorphic Jacobi forms of weight $0$ and weakly holomorphic vector-valued modular forms of weight $-\mathrm{rank}(K) / 2$ for the Weil representation for $K$. The latter may be interpreted as input functions in Borcherds' sense after identifying the Weil representations for $K$ and $L = K \oplus 2U$.

\section{Proof of Theorem 1}

(i) $\Rightarrow$ (ii): Suppose $\Psi$ is the Borcherds lift of a weakly holomorphic Jacobi form $$\phi(\tau, z) = \sum_{n \gg -\infty} \sum_{\ell \in K'} c(n, \ell) q^n \zeta^{\ell},$$ and let $N$ be the greatest common divisor of all inner products of vectors in $K$. Borcherds' congruence (\cite{Borcherds}, Theorem 11.2; a simpler proof for Jacobi forms is given in Corollary 4.5 of \cite{Wang}) implies $$N \sum_{\ell \in K'} c(0, \ell) \equiv 0 \, (\text{mod} \, 24).$$ By definition, each term $c(0, \ell)$ in the above sum is integral. If $\Psi$ has half-integral weight, then $c(0,0)$ is odd. Since $c(0, \ell)$ is symmetric under $\ell \mapsto -\ell$, it follows that $\sum_{\ell \in K'} c(0, \ell)$ is odd and therefore $8 | N$. \\

(ii) $\Rightarrow$ (i): The existence proof is constructive. We proceed in two steps. \\

Step 1. We first construct a  Borcherds product of half-integral weight when $K = 8\mathbb{Z}^N$; in other words, a Gram matrix for $K$ is given by the diagonal matrix $\mathrm{diag}(8,...,8)$. Recall (e.g. \cite{Borcherds1}, Example 5 of section 15; or \cite{GN}, Example 2.3) that there is a holomorphic product of weight $1/2$ attached to the lattice $K = 8 \mathbb{Z}$. The input function as a weak Jacobi form can be given as a quotient of theta functions: 
$$
\phi_{0,4}(\tau, z) = \frac{\vartheta(\tau, 3z)}{\vartheta(\tau, z)} = \zeta + 1 + \zeta^{-1} + O(q),
$$ 
where $q = e^{2\pi i \tau}$, $\zeta = e^{2\pi i z}$, and (by Jacobi's triple product)
\begin{align*} \vartheta(\tau, z) &= \sum_{n = -\infty}^{\infty} \left( \frac{-4}{n} \right) q^{\frac{n^2}{8}} \zeta^{\frac{n}{2}} \\ &=q^{\frac{1}{8}}(\zeta^{\frac{1}{2}}-\zeta^{-\frac{1}{2}})\prod_{n=1}^\infty(1-q^n\zeta)(1-q^n\zeta^{-1})(1-q^n). \end{align*}
For any $N \in \mathbb{N}$, we obtain a weak Jacobi form of weight $0$ and index $8 \mathbb{Z}^N$ by repeatedly taking the direct product of $\phi_{0,4}$ with itself: $$\phi_N(\tau, z_1,...,z_n) := \phi_{0,4}(\tau, z_1) \phi_{0,4}(\tau, z_2) ... \phi_{0,4}(\tau, z_n).$$ The Borcherds lift $\Psi_N$ of $\phi_N$ is then a (meromorphic) Borcherds product of weight $1/2$. By a theorem of Bruinier (\cite{Bruinier}, Theorem 1.3) $\Psi_N$ is a quotient of two holomorphic Borcherds products; one of these is of half-integral weight. \\

Step 2. Let $L$ be any even lattice satisfying (ii).  The rescaled lattice $K(1/8)$ is integral (not necessarily even) and embeds in a positive-definite Type I unimodular lattice of rank at most $l + 1$ by Corollary 8 of \cite{CS}. It follows that $K$ is embedded in an even lattice $M$ lying in the same genus as $8\mathbb{Z}^N$ for some $N$. Since even lattices in the same genus yield equivalent Jacobi forms, we obtain from Step 1 a holomorphic half-integral weight Borcherds product attached to $M \oplus 2U$. By \cite{Ma}, its quasi-pullback to $L$ is a Borcherds product of half-integral weight.

\section{Examples}

Every step of the above proof can be made constructive, but the half-integral weight products constructed in this way tend to have rather large weight and complicated divisors. Therefore it seems worthwhile to present a few simpler examples of half-integral weight Borcherds products on lattices of rank greater than $5$.\\

In this section we give two examples of half-integral weight Borcherds products associated to lattices of signature $(4, 2)$. These can also be interpreted as Hermitian modular forms of degree two over the Eisenstein and Gaussian integers, respectively. In both cases we give only the principal part of the input form as a vector-valued modular form (which determines the weight and divisor). The two examples were computed in SAGE \cite{sagemath}.

\begin{ex} Take $K = \mathbb{Z}^2$ with Gram matrix $\begin{psmallmatrix} 16 & 8 \\ 8 & 16 \end{psmallmatrix}$. The lattice $L = K \oplus 2 U$ admits (at least) four holomorphic Borcherds products of weight $9/2$. These four forms are conjugates under the symmetries of the discriminant form. The principal part of one of the forms is 
\begin{align*} 
9 \mathfrak{e}_{(0,0)} &+ q^{-1/24} (4 \mathfrak{e}_{\pm (11/12, 1/24)} + 4 \mathfrak{e}_{\pm (1/24, 1/24)} + 4 \mathfrak{e}_{\pm (1/24, 11/12)}  - \mathfrak{e}_{\pm (7/12, 5/24)} \\
&- \mathfrak{e}_{\pm (5/24, 5/24)} - \mathfrak{e}_{\pm (5/24, 7/12)}  + 2 \mathfrak{e}_{\pm (5/12, 7/24)} + 2 \mathfrak{e}_{\pm (7/24, 7/24)} + 2 \mathfrak{e}_{\pm (7/24, 5/12)}\\
& + \mathfrak{e}_{\pm (1/12, 11/24)} + \mathfrak{e}_{\pm (11/24, 11/24)} + \mathfrak{e}_{\pm (11/24, 1/12)})  + q^{-1/8} (3 \mathfrak{e}_{\pm (1/8, 0)}\\
& + 3 \mathfrak{e}_{\pm (0, 1/8)}+ 3 \mathfrak{e}_{\pm (1/8, 7/8)}) + q^{-1/6} (\mathfrak{e}_{\pm (1/6, 5/12)} + \mathfrak{e}_{\pm (5/12, 5/12)} \\
&+ \mathfrak{e}_{\pm (5/12, 1/6)} + 2 \mathfrak{e}_{\pm (5/6, 1/12)} +  2\mathfrak{e}_{\pm (1/12, 1/12)} + 2 \mathfrak{e}_{\pm (1/12, 5/6)}), 
\end{align*} 
where $\mathfrak{e}_{\pm \gamma}$ represents $\mathfrak{e}_{\gamma} + \mathfrak{e}_{-\gamma}$.
\end{ex}

\begin{ex} Take $K = \mathbb{Z}^2$ with Gram matrix $\begin{psmallmatrix} 8 & 0 \\ 0 & 8 \end{psmallmatrix}.$ The lattice $L = K \oplus 2U$ admits two holomorphic Borcherds products of weight $7/2$. These can be interpreted as Hermitian modular forms for a level $4$ subgroup over the Gaussian integers. These forms are conjugate under an involution of the discriminant form. The principal part of one of the forms is \begin{align*} 7 \mathfrak{e}_{(0, 0)} &+ q^{-1/16} ( 3\mathfrak{e}_{\pm (1/8, 0)} + 3\mathfrak{e}_{\pm (0, 1/8)}  + \mathfrak{e}_{\pm (1/2, 1/8)} +  \mathfrak{e}_{\pm (1/8, 1/2)}) \\ &+ q^{-1/8} (\mathfrak{e}_{\pm (1/8, 1/8)} + \mathfrak{e}_{\pm (1/8, 7/8)}) \\ &+q^{-1/4} (\mathfrak{e}_{\pm (1/4, 0)} + \mathfrak{e}_{\pm (0, 1/4)}). \end{align*}
\end{ex}

\bibliographystyle{plainnat}
\bibliofont
\bibliography{\jobname}

\end{document}